\documentclass[10pt,twoside]{article}
\usepackage{times, fancyhdr}
\usepackage{cite}
\usepackage[ ]{graphicx}
\usepackage{epsfig}
\usepackage{amsmath}
\usepackage{amssymb}
\usepackage{amsfonts}
\usepackage{latexsym,amsmath,amsfonts,amssymb,amsthm}

\usepackage{amsthm}
\usepackage{indentfirst}
\usepackage{sectsty}
\sectionfont{\large}
\subsectionfont{\normalsize\bf}

\usepackage{float}

 \numberwithin{equation}{section}

\newtheorem{theorem}{Theorem}[section]
\newtheorem{lemma}[theorem]{Lemma}

\title{\vspace{-0.844in}\parbox{\linewidth}{\footnotesize\noindent
{International Journal of Modern Mathematics, 3(2008), no.\,2,
197--206.} \hspace{\stretch{0}} }
{\large {\bf On Sun's Conjecture concerning Disjoint Cosets}}}
\date{}

\author{\normalsize Wan-Jie Zhu}
\begin{document}

\maketitle

\thispagestyle{empty}

\renewcommand{\theequation}{\thesection.\arabic{equation}}
\makeatletter \@addtoreset{equation}{section}

\date{}{}
\def\colon{{:}\;}
\def\Z{\Bbb Z}
\def\zp{\Z^+}
\def\Q{\Bbb Q}
\def\N{\Bbb N}
\def\R{\Bbb R}
\def\C{\Bbb C}
\def\al{\alpha}
\def\l{\left}
\def\r{\right}
\def\bg{\bigg}
\def\({\bg(}
\def\){\bg)}
\def\[{\bg\lfloor}
\def\]{\bg\rfloor}
\def\t{\mbox}
\def\f{\frac}
\def\si{\sim_2}
\def\pr{\prec_2}
\def\pe{\preceq_2}
\def\B{B_{k+1}}
\def\E{E_k}
\def\P{\Z_q}
\def\q{\ (\t{\rm mod}\ q)}
\def\qc{q\t{-adic}}
\def\se {\subseteq}
\def\sp {\supseteq}
\def\sm{\setminus}
\def\Ar{\Longrightarrow}
\def\La{\Longleftrightarrow}
\def\bi{\binom}
\def\eq{\equiv}
\def\cs{\ldots}
\def\ls{\leqslant}
\def\gs{\geqslant}
\def\mo{\t{\rm mod}}
\def\Tor{\t{\rm Tor}}
\def\sign{\t{\rm sign}}
\def\o{\t{\rm ord}}
\def\ch{\t{\rm ch}}
\def\per{\t{\rm per}}
\def\ve{\varepsilon}
\def\da{\delta}
\def\Da{\Delta}
\def\la{\lambda}
\def\ta{\theta}
\def\ord{\t{\rm ord}}
\def\Tor{\t{\rm Tor}}
\def\kl{k\varphi(q^{\al})+r}
\def\mx{\langle-x\rangle_q}
\def\Mr{\sum\Sb 0\ls j<q\\m\mid j-r-qs\endSb(x+j)^k}
\def\qed{}
\def\Proof{\noindent{\it Proof}}
\def\Def{\medskip\noindent{\it Definition}\ }
\def\Remark{\medskip\noindent{\it  Remark}}
\def\Ack{\medskip\noindent {\bf Acknowledgment}}

\begin{abstract}
\noindent In 2004, Zhi-Wei Sun posed the following conjecture: If
$a_1G_1,\ldots,a_kG_k\ (k>1)$ are finitely many pairwise disjoint
left cosets in a group $G$ with all the indices $[G:G_i]$ finite,
then for some $1\le i<j\le k$, the greatest common divisor of
$[G:G_i]$ and $[G:G_j]$ is at least $k$. In this paper, we confirm
Sun's conjecture for $k=3,4$.


\end{abstract}

\section{Introduction}

Let $H$ be any subgroup of a (multiplicative) group $G$. A left
coset of $H$ has the form $aH=\{ah:\ h\in H\}$ with $a\in G$.
$[G:H]$, the index of $H$ in $G$, is the cardinality of the set
$G/H=\{aH:\ a\in G\}$. If $k=[G:H]<\infty$, then we can partition
$G$ into $k$ distinct left cosets of $H$ in $G$.

In 2004, Zhi-Wei Sun proposed the following conjecture on disjoint
cosets.
\medskip

 \noindent {\bf Sun's Conjecture} (\cite [Conjecture 1.2] {S06}). {\it  Let
$a_1G_1,\ldots,a_kG_k\ (k>1)$ be finitely many pairwise disjoint
left cosets in a group $G$ with $[G:G_i]<\infty$ for all
$i=1,\ldots,k$. Then, for some $1\le i<j\le k$, the greatest
common divisor $([G:G_i],[G:G_j])$ of $[G:G_i]$ and $[G:G_j]$ is
at least $k$. }
\medskip

By \cite [Remark 2.2]{S06} and \cite [Remark 1.5(b)]{S06}, Sun's
conjecture holds when $k=2$ or $G$ is a $p$-group with $p$ a
prime. K. O'Bryant \cite{O} proved Sun's conjecture in the special case $G=\mathbb Z$ and $k\le 20$.

 In this paper we confirm Sun's conjecture for $k\in\{3,4\}$.

 \begin{theorem} Let
$a_1G_1,\ldots,a_kG_k\ (k\in\{3,4\})$ be pairwise disjoint
left cosets in a group $G$ with $[G:G_i]<\infty$ for all
$i=1,\ldots,k$. Then
$([G:G_i],[G:G_j])\ge k$ for some $1\le i<j\le k$.
\end{theorem}

\section{The case $k=3$}

\begin{lemma} Let $H$ and $K$ be two subgroups of a group $G$.

{\rm (i) (\cite[p. 41]{R})} $HK$ is a subgroup of $G$ if and only if
$HK=KH$.

{\rm (ii) (\cite[Chapter 1]{Suz})} $HK$ contains exactly $[H:H\cap K]$ left cosets of $K$.

{\rm (iii) (\cite[Lemma 2.1]{S06})} If $[G:H]$ and $[G:K]$ are
finite and relatively prime, then $HK$ coincides with $G$.

{\rm (iv) (\cite[Lemma 2.1(i)]{S01})} $HK=G$ if and only if $xH\cap
yK\not=\emptyset$ for all $x,y\in G$.

{\rm (v)} Suppose that $HK=KH$ and $xH\cap yK=\emptyset$, where
$x,y\in G$. Then  $xHK\cap yHK=\emptyset$.
\end{lemma}

\Proof. Parts (i)-(iv) are known. So we just prove part (v). Assume
that $xHK\cap yHK\not=\emptyset$. Then $x^{-1}y\in HK$ and hence
$x^{-1}y=hk$ for some $h\in H$ and $k\in K$. It follows that
$xh=yk^{-1}\in xH\cap yK$, which contradicts the condition $xH\cap
yK=\emptyset$. We are done.

\begin{theorem}Let $a_1G_1,a_2G_2$ and $a_3G_3$ be pairwise
disjoint left cosets in a group $G$ with $[G:G_i]<\infty$ for
$i=1,2,3$. Then $([G:G_i],[G:G_j])\ge 3$ for some $1\le i<j\le 3$.
\end{theorem}

\Proof. Suppose that we don't have the desired result. Then,
whenever $1\le i<j\le3$, we have $([G:G_i],[G:G_j])\le 2$ and
hence $([G:G_i],[G:G_j])=2$ (otherwise $a_iG_i\cap
a_jG_j\not=\emptyset$ by Lemma 2.1(ii)-(iv)).

\smallskip
Write $[G:G_i]=2q_i$ with $q_i\in \Z^+=\{1,2,\ldots\}$. Fix $1\le
i\not=j \le3$. As $([G:G_i],[G:G_j])=2$, $q_i$ is relatively prime
to $q_j$. Since $a_iG_i\cap a_jG_j=\emptyset$, $G_iG_j\not=G$ by
Lemma 2.1(iv), and hence $[G_i:G_i\cap G_j]<[G:G_j]=2q_j$  by
Lemma 2.1(ii). As both $2q_i=[G:G_i]$ and $2q_j=[G:G_j]$ divide
$[G:G_i\cap G_j]$, $2q_iq_j$ divides $[G:G_i\cap
G_j]=[G:G_i][G_i:G_i\cap G_j]$ and hence $q_j\mid [G_i:G_i\cap
G_j]$. Therefore $[G_i:G_i\cap G_j]=q_j$.

Let $k\in \{1,2,3\}\backslash\{i,j\}$. As the above, we also have
$[G_i:G_i\cap G_k]=q_k$. Since
$$([G_i:G_i\cap G_j],[G_i:G_i\cap
G_k])=(q_j,q_k)=1,$$
 applying Lemma 2.1(iii) we find that
$$G_i=(G_i\cap G_j)(G_i\cap G_k)=(G_i\cap G_k)(G_i\cap G_j).$$
Similarly, $$G_j=(G_i\cap G_j)(G_j\cap G_k)=(G_j\cap G_k)(G_i\cap
G_j)$$ and $$G_k=(G_i\cap G_k)(G_j\cap G_k)=(G_j\cap G_k)(G_i\cap
G_k).$$
Thus
$$\aligned G_iG_j=&(G_i\cap G_k)(G_i\cap G_j)(G_i\cap G_j)(G_j\cap G_k)\\
=&(G_i\cap G_k)(G_i\cap G_j)(G_j\cap G_k).\endaligned$$

\noindent Since any two of the subgroups $G_i\cap G_j,\ G_i\cap
G_k,\ G_j\cap G_k$ are commutable, $G_iG_j$ coincides with
$H=(G_1\cap G_2)(G_1\cap G_3)(G_2\cap G_3)$, which is a subgroup of
$G$.

As $a_iG_i\cap a_jG_j=\emptyset$, by Lemma 2.1(v) we have
$a_iG_iG_j\cap a_jG_iG_j=\emptyset$, i.e., $a_iH\cap a_jH=\emptyset$.
Note that
$$[G:H]=[G:G_iG_j]=\frac{[G:G_j]}{[G_iG_j:G_j]}=\frac{[G:G_j]}{[G_i:G_i\cap G_j]}=\frac{2q_j}{q_j}=2.$$
On the other hand, $a_1H,a_2H$ and $a_3H$ are pairwise
disjoint. So we get a contradiction.

\section{The case $k=4$}

\noindent In this section, we prove the following result.

\begin{theorem} Let $a_1G_1,a_2G_2,a_3G_3$ and $a_4G_4$ be left cosets in a group $G$
with $[G:G_i]<\infty$ for $i=1,2,3,4$. Then $([G:G_i],[G:G_j])\ge 4$ for some $1\le i<j\le 4$.
\end{theorem}

\begin{lemma} Let $H_1,\ldots,H_k$ be $k$ subgroups of a group
$G$. Then there is a bijection from
$$S=\{(C_1,\ldots,C_k)\in
G/H_1\times\cdots\times G/H_k:\ C_1\cap\cdots\cap
C_k\not=\emptyset\}$$ to $G/\bigcap_{i=1}^kH_i.$
\end{lemma}

\Proof. For $(C_1,\ldots,C_k)\in S$. Let
$\sigma(C_1,\ldots,C_k)=\bigcap_{i=1}^kC_i$. For $x\in
\bigcap_{i=1}^kC_i$, we have $C_i=xH_i$ for $i=1,\ldots,k$, and
hence $\bigcap_{i=1}^kC_i=\bigcap_{i=1}^kxH_i=xH$ where $H=\bigcap
_{i=1}^kH_i$. Note also that $xH=\bigcap_{i=1}^kxH_i=\sigma
(xH_1,\ldots,xH_k)$ for any $x\in G$. So $\sigma$ is surjective.
 If $(x_1H_1,\ldots,x_kH_k)$, $(y_1H_1,\ldots,y_kH_k)\in S$ and
 $\bigcap_{i=1}^kx_iH_i=\bigcap_{i=1}^ky_iH_i$, then for
 $z\in \bigcap_{i=1}^kx_iH_i=\bigcap_{i=1}^ky_iH_i$, we have
 $x_iH_i=zH_i=y_iH_i$
  for $i=1,\ldots,k$. Thus $\sigma$ is also injective and hence it
  is a bijection form $S$ to $G/H$.

\Remark\  3.1. Let $G_1$ and $G_2$ be subgroups of a group $G$, and
let $H_1$ and $H_2$ be subgroups of $G_1$ and $G_2$ with finite
index, respectively, satisfying $H_1\cap H_2=G_1\cap G_2$. In view
of Lemma 3.2, there exists a bijection
$$\sigma:  S=\{(C_1,C_2)\in
G/G_1\times G/G_2:\ C_1\cap C_2\not=\emptyset\}\rightarrow G/(H_1\cap
H_2)$$
Suppose $aH_1$ and $\tilde aH_1$ are distinct left cosets of $H_1$ contained in
a left coset of $G_1$. Then
$$\sigma (aG_1,\tilde aG_2)=aG_1\cap \tilde aG_2=\tilde aG_1\cap \tilde aG_2=\tilde a(H_1\cap H_2).$$
Thus
$$aH_1\cap\tilde aH_2\subseteq aG_1\cap\tilde aG_2=\tilde a(G_1\cap G_2)=\tilde a (H_1\cap H_2)\se \tilde aH_1.$$
But $aH_1\cap \tilde aH_1=\emptyset$, so $aH_1\cap\tilde aH_2=\emptyset$.

\medskip
\begin{lemma} Let $H$ and $K$ be subgroups of a group $G$.
Then $$[K:H\cap K]=|\{C\in G/H:\ K\cap C\not=\emptyset\}|.$$
\end{lemma}

\Proof. Let $g\in G$. If $gh=k$ with $h\in H$ and $k\in K$, then
$gH=kh^{-1}H=kH\subseteq KH$. If $gH\subseteq KH$, then we can write
 $g=kh^{-1}$ with $h\in H$ and $k\in K$, and hence $gh=k\in gH\cap K$.
 So $gH\cap K\not=\emptyset$ if and only if $gH\subseteq KH$.

In view of Lemma 2.1(ii) and the above, we have $$[K:H\cap K]=|\{gH:\ g\in G\ \&\ gH\subseteq
KH\}|=|\{C\in G/H:\ C\cap K\not=\emptyset\}|.$$ This concludes the
proof.

\bigskip
\noindent {\it Proof of Theorem 3.1}. Suppose that the desired
result is false. We want to deduce a contradiction. For any
$i,j\in\{1,2,3,4\}$ with $i\not=j$, clearly
$d_{ij}=([G:G_i],[G:G_j])\in\{2,3\}$, since $a_iG_i\cap
a_jG_j=\emptyset$. Let $i,j$ and $k$ be three distinct elements of
$\{1,2,3,4\}$. As $d_{ij},d_{ik}$ and $d_{jk}\in\{2,3\}$, two of
$d_{ij},d_{ik},d_{jk}$ are equal, say,
$d_{ij}=d_{ik}=d\in\{2,3\}$. As $d$ divides both $[G:G_j]$ and
$[G:G_k]$, we have $d\mid d_{jk}$. Note that $d_{jk}\le 3<2d$ and
hence $d_{jk}=d$. Since $a_iG_i,a_jG_j$ and $a_kG_k$ are pairwise
disjoint, $\max\{d_{ij},d_{ik},d_{jk}\}\ge 3$ by Theorem 2.1, and
thus $d_{ij}=d_{ik}=d_{jk}=d=3$.

Write $[G:G_i]=3q_i$ for $i=1,2,3,4$. Then $q_1,q_2,q_3,q_4$ are pairwise coprime. Let $1\le i<j\le 4$.
As $[G:G_i\cap G_j]$ is divisible by $[G:G_i]=3q_i$ and $[G:G_j]=3q_j$,
we may write $[G:G_i\cap G_j]=3q_iq_jr_{ij}$ with $r_{ij}\in\Z^+$.
As $a_iG_i\cap a_jG_j=\emptyset$, we have $G_iG_j\not=G$ by Lemma 2.1(iv) and hence
$$[G_i:G_i\cap G_j]=|\{gG_j:\, g\in G\ \&\ gG_j\subseteq G_iG_j\}|<[G:G_j]$$
by Lemma 2.1(ii). Thus $[G:G_i\cap G_j]<[G:G_i][G:G_j]$. So $3q_iq_jr_{ij}<3q_i3q_j$
and hence $r_{ij}\le 2$.

Let $1\le i<j<k\le 4$. We may write $[G:G_i\cap G_j\cap
G_k]=3q_iq_jq_kr_{ijk}$ with $r_{ijk}\in \Z^+$.
Observe that
$$\aligned &[G_i\cap G_j:G_i\cap G_j\cap G_k]\\
=&|\{gG_k:\ g\in G\ \&\ gG_k\subseteq (G_i\cap G_j)G_k\}|\\
\le&|\{gG_k:\ g\in G\ \&\ gG_k \subseteq G_iG_k\}|\\
=&[G_i:G_i\cap G_k]\endaligned$$
and hence
$$\frac{[G:G_i\cap G_j\cap G_k]}{[G:G_i\cap G_j]}\le \frac{[G:G_i\cap G_k]}{[G:G_i]}.$$
In other words,
$$\frac{3q_iq_jq_kr_{ijk}}{3q_iq_jr_{ij}}\le\frac{3q_iq_kr_{ik}}{3q_i},$$
i.e.,
\begin{equation}\label{eq-3.1}r_{ijk}\le r_{ij}r_{ik}. \end{equation}
It follows that $r_{ijk}\le r_{ij}r_{ik}\le 2\times 2=4$.
(\ref{eq-3.1}) is useful and it was suggested by Prof. Zhi-Wei
Sun.

Let $1\le i<j<k\le 4$. Set $S=G/G_i\times G/G_j\times G/G_k$, and
$$\aligned &S_{ij}=\{(C_i,C_j,C_k)\in S:\ C_i\cap C_j\not=\emptyset\},
\\&S_{ik}=\{(C_i,C_j,C_k)\in S:\ C_i\cap C_k\not=\emptyset\},
\\&S_{jk}=\{(C_i,C_j,C_k)\in S:\ C_j\cap C_k\not=\emptyset\}.\endaligned$$
Then $|S|=[G:G_i][G:G_j][G:G_k]=3q_i3q_j3q_k=27q_iq_jq_k$.

\noindent Also,
$$\aligned |S_{ij}|=&|\{(C_i,C_j)\in G/G_i\times G/G_j:\ C_i\cap C_j\not=\emptyset\}|\times|G/G_k|\\
=&[G:G_i\cap G_j][G:G_k]=9q_iq_jq_kr_{ij}.\endaligned$$
Similarly,  $|S_{ik}|=9q_iq_jq_kr_{ik}$ and
$|S_{jk}|=9q_iq_jq_kr_{jk}$.

Clearly $$|S_{ij}\cap S_{ik}|=|\{(C_i,C_j,C_k)\in S:\
C_i\cap C_j\not=\emptyset \ \&\ C_i\cap C_k\not=\emptyset\}|$$

\noindent Given $C_i\in G/G_i$ and $C_j\in G/G_j$ with $C_i\cap C_j$
containing some $x\in G$, we have
$$\aligned
&|\{C_k\in G/G_k:\ C_k\cap C_i\not=\emptyset\}|\\
=&|\{xC_k:\ C_k\in G/G_k\ \&\ xC_k\cap C_i\not=\emptyset\}|\\
=&|\{xC_k:\ C_k\in G/G_k\ \&\ xC_k\cap xG_i\not=\emptyset\}|\\
=&|\{xC_k:\ C_k\in G/G_k\ \&\ C_k\cap G_i\not=\emptyset\}|\\
=&|\{C_k\in G/G_k:\ C_k\cap G_i\not=\emptyset\}|\\
=&[G_i:G_i\cap G_k]\ \ \mbox{(by Lemma 3.3)}.
\endaligned$$
Thus
$$\aligned |S_{ij}\cap S_{ik}|=&|\{(C_i,C_j)\in G/G_i\times G/G_j: C_i\cap C_j\not=\emptyset\}|[G_i:G_i\cap G_k]
\\=&\frac{[G:G_i\cap G_j][G:G_i\cap G_k]}{[G:G_i]}=3q_iq_jq_kr_{ij}r_{ik}.\endaligned$$
Similarly,  $|S_{ij}\cap S_{jk}|=3q_iq_jq_kr_{ij}r_{jk}$
and $|S_{ik}\cap S_{jk}|=3q_iq_jq_kr_{ik}r_{jk}$.

Observe that $$\aligned &|S_{ij}\cap S_{ik}\cap
S_{jk}|\\=&|\{(C_i,C_j,C_k)\in S: \ C_i\cap C_j\not=\emptyset,\
C_i\cap C_k\not=\emptyset,\ C_j\cap C_k\not=\emptyset\}|\\\ge
&|\{(C_i,C_j,C_k)\in S:\ C_i\cap C_j\cap
C_k\not=\emptyset\}|\\=&[G:G_i\cap G_j\cap
G_k]=3q_iq_jq_kr_{ijk}\ \ \mbox{(by Lemma 3.2)}.\endaligned$$

 By the inclusion-exclusion principle,
$$\aligned N_{ijk}=&|\{(C_i,C_j,C_k)\in S:\ C_i,C_j,C_k\ \text{are pairwise disjoint}\}|
\\=&|\{(C_i,C_j,C_k)\in S:\ (C_i,C_j,C_k)\notin S_{ij}\cup S_{ik}\cup S_{jk}\}|
\\=&|S|-|S_{ij}\cup S_{jk}\cup S_{jk}|
\\=&|S|-(|S_{ij}|+|S_{ik}|+|S_{jk}|)+(|S_{ij}\cap S_{ik}|+|S_{ij}\cap S_{jk}|
+|S_{ik}\cap S_{jk}|)\\&-|S_{ij}\cap S_{ik}\cap S_{jk}|
\\\le &27q_iq_jq_k-9q_iq_jq_k(r_{ij}+r_{ik}+r_{jk})+3q_iq_jq_k(r_{ij}r_{ik}+r_{ij}r_{jk}+r_{ik}r_{jk})
\\&-3q_iq_jq_kr_{ijk}.\endaligned$$
Since $a_iG_i,a_jG_j,a_kG_k$ are pairwise disjoint,
$N_{ijk}>0$ and hence
$$27q_iq_jq_k-9q_iq_jq_k(r_{ij}+r_{ik}+r_{jk})+3q_iq_jq_k(r_{ij}r_{ik}+r_{ij}r_{jk}+r_{ik}r_{jk})
-3q_iq_jq_kr_{ijk}>0.$$
Thus
\begin{equation}\label{eq-3.2}r_{ijk}<9-3(r_{ij}+r_{ik}+r_{jk})+(r_{ij}r_{ik}+r_{ij}r_{jk}+r_{ik}r_{jk}).\end{equation}

When $\{r_{ij},r_{ik},r_{jk}\}=\{1,2,r\}$, (3.2) gives
$$r_{ijk}<9-3(1+2+r)+(1\times 2+1\times r+2\times r)=2.$$
So we have
\begin{equation}\label{eq-3.3}\{r_{ij},r_{ik},r_{jk}\}=\{1,2,r\}\Rightarrow r_{ijk}=1.\end{equation}
If $r_{ij}=r_{ik}=r_{jk}=r$, then (3.2) yields that
$$r_{ijk}<9-3\times 3r+3r^2=3(r^2-3r+3)=3((r-1)(r-2)+1)=3.$$
So we always have $r_{ijk}\le 2$.

Since $[G:G_i\cap G_j]\mid [G:G_i\cap G_j\cap G_k]$,
$3q_iq_jr_{ij}$ divides $3q_iq_jq_kr_{ijk}$ and hence
\begin{equation}\label{eq-3.4}r_{ij}\mid q_kr_{ijk}.\end{equation}
If $2\nmid q_k$, then $r_{ij}\mid r_{ijk}$,  since
$(q_k,r_{ij})=1$.
Suppose $q_i,q_j,q_k$ are all odd, we then have
$$r_{ij}=r_{ik}=r_{jk}=r_{ijk}.$$
In fact,$$\aligned {\rm (a)}\ &\text{if\ }
r_{ij}=r_{ik}=r_{jk}=1,\text{ then }r_{ijk}\le r_{ij}r_{ik}=1;
\\{\rm (b)}&\text{ if }r_{ij}=r_{ik}=r_{jk}=2, \mbox{then}\ r_{ijk}=2\text{, since }r_{ijk}\le 2
\text{ and }r_{ij}\mid r_{ijk};\\{\rm (c)}&\text{ if }
\{r_{ij},r_{ik},r_{jk}\}=\{1,2,r\},\text{ then }r_{ijk}=1\ \mbox{by (3.3) and also}\
2\mid r_{ijk}\text{by the above.}\endaligned$$

 Since $q_1,q_2,q_3,q_4$ are pairwise coprime, there are no
two even numbers among $q_1,q_2,q_3,q_4$. Without loss of
generality, we assume that $q_1,q_2,q_3$ are all odd. By the
above, we have $r_{12}=r_{13}=r_{23}=r_{123}\in\{1,2\}$.

\medskip
{\it Case} 1. $r_{12}=r_{13}=r_{23}=r_{123}=1$.

Let $i,j,k$ be any permutation of $1,2,3$. Note that
$$[G_i:G_i\cap G_j]=\frac{[G:G_i\cap G_j]}{[G:G_i]}=\frac{3q_iq_jr_{ij}}{3q_i}=
q_jr_{ij}=q_j.$$ Similarly,
$[G_i:G_i\cap G_k]= q_kr_{ik}=q_k$ and $[G_i:G_i\cap
G_4]=q_4r_{i4}$. Since any two of $q_j,q_k,q_4r_{i4}$ are coprime,
no matter $r_{i4}$ is $1$ or $2$, by Lemma 2.1(iii) we have
\begin{equation}\label{eq-3.5}G_i=(G_i\cap G_j)(G_i\cap G_k)=(G_i\cap G_j)(G_i\cap
G_4)=(G_i\cap G_k)(G_i\cap G_4).\end{equation}

 In view of (3.5),
$$G_1G_2=G_2G_1=G_1G_3=G_3G_1=G_2G_3=G_3G_2=(G_1\cap G_2)(G_1\cap G_3)(G_2\cap G_3).$$
Denote this subgroup by $H=G_iG_j\ (1\le i<j\le 3)$. Then
$$[G:H]=[G:G_1G_2]=\frac{[G:G_2]}{[G_1:G_1\cap G_2]}=\frac{3q_2}{q_2}=3.$$
As $a_1G_1,a_2G_2,a_3G_3$ are pairwise disjoint, so are
$a_1H,a_2H,a_3H$ by Lemma 2.1(v).

 Suppose that $a_iH\cap a_4H\not=\emptyset$ for some $i=1,2,3$. Then $a_4^{-1}a_i\in H$.
 Take $j\in\{1,2,3\}\sm\{i\}$ and note that
 $$H=G_jG_i=(G_j\cap G_4)(G_j\cap G_i)G_i=(G_j\cap G_4)G_i.$$
 So there are $ g_i\in G_i$ and $ g_4\in
G_j\cap G_4$ such that $a_4^{-1}a_i=g_4g_i$ and hence
$$a_ig_i^{-1}=a_4g_4\in a_iG_i\cap a_4G_4,$$
contradicting the known
condition $a_iG_i\cap a_4G_4=\emptyset$.

By the above discussion, $a_1H,a_2H,a_3H,a_4H$ are
distinct left cosets of $H$, which contradicts $[G:H]=3$.

\medskip {\it Case} 2. $r_{12}=r_{13}=r_{23}=r_{123}=2$.

Suppose that $r_{14},r_{24},r_{34}$ are not all equal, say,
$r_{14}=1$ and $r_{24}=2$. Then $q_4$ must be even, otherwise we
will get $1=r_{14}=r_{12}=r_{24}=2$ since $q_1, q_2,q_4$ are odd.
We also have $r_{124}=1$ by (\ref{eq-3.2}). With the help of
(3.4), $r_{24}=2$ divides $q_1r_{124}=q_1$. So we get a
contradiction.

Let $1\le i<j\le 3$. Recall that $r_{ij4}\le 2$.
When $r_{14}=r_{24}=r_{34}=1$,  $q_4$ is even (otherwise $q_1,q_2,q_4$ are odd and hence $r_{14}=r_{24}=r_{12}=2$)
 and also $r_{ij4}\le
r_{i4}r_{j4}=1$.  By (3.4), $r_{j4}$ divides
$q_ir_{ij4}$. Since $2\nmid q_i$, if $r_{14}=r_{24}=r_{34}=2$ then $r_{j4}=2=r_{ij4}$.

By the above, $r_{14}=r_{24}=r_{34}=r_{ij4}$ for any
$1\le i<j\le 3$.

Fix $1\le i<j\le 4$ and write $\{1,2,3,4\}\sm\{i,j\}=\{k,l\}$.
If $j\le 3$ then
$$[G_i\cap G_j:G_i\cap G_j\cap G_k]=\frac{[G:G_i\cap G_j\cap G_k]}{[G:G_i\cap G_j]}$$
divides
$$\f{3q_iq_jq_k\times2}{3q_iq_jr_{ij}}=q_k;$$
if $j=4$ then we have
$$[G_i\cap G_j:G_i\cap G_j\cap G_k]=\f{3q_iq_kq_4r_{ik4}}{3q_iq_4r_{i4}}=q_k.$$
Since $q_k$ and $q_l$ are coprime, so are $[G_i\cap G_j:G_i\cap G_j\cap G_k]$ and $[G_i\cap G_j:G_i\cap G_j\cap G_l]$.
Thus, by Lemma 2.1(iii) we have
$$G_i\cap G_j=(G_i\cap G_j\cap G_k)(G_i\cap G_j\cap G_l).$$

Let $1\le k \le4$. For $i,j\in\{1,2,3,4\}\sm\{k\}$ with $i\not=j$, it follows from the above that
\begin{equation}\label{eq-3.6}{(G_k\cap G_i)(G_k\cap G_j)}
=\prod_{\{s,t\}\subseteq\{1,2,3,4\}\backslash \{k\},\ s<t}{(G_k\cap
G_s\cap G_t)}\end{equation}
Denote by $H_k$ the subgroup of $G_k$ given by (\ref{eq-3.6}).
Choose $i,j\in\{1,2,3\}\sm\{k\}$ with $i\not=j$. Then
\begin{align*}[G:H_k]=&\f{[G:G_k\cap G_j]}{[(G_k\cap G_i)(G_k\cap G_j):G_k\cap G_j]}=\f{[G:G_k\cap G_j]}{[G_k\cap G_i:G_i\cap G_j\cap G_k]}
\\=&\f{[G:G_i\cap G_k][G:G_j\cap G_k]}{[G:G_i\cap G_j\cap G_k]}=\f{3q_iq_kr_{ik}3q_jq_kr_{jk}}{3q_iq_jq_kr_{ijk}}
=3q_k\f{r_{ik}r_{jk}}{r_{ijk}}
\\=& \begin{cases} 3q_k&\text{if}\ k=4\ \&\ r_{14}=r_{24}=r_{34}=1,
\\6q_k&\text{otherwise}.
\end{cases}\end{align*}
Thus, when $k\not=4$ we have
$$[G_k:H_k]=\f{[G:H_k]}{[G:G_k]}=\f{6q_k}{3q_k}=2$$
and hence $a_kG_k\sm a_kH_k=\tilde a_kH_k$ for some $\tilde a_k\in G$.

Let $1\le i\not=j\le4$. It is easy to see that
\begin{equation}\label{eq-3.7}H_iH_j=(G_1\cap
G_2\cap G_3)(G_1\cap G_2\cap G_4)(G_1\cap G_3\cap G_4)(G_2\cap
G_3\cap G_4).\end{equation}
We denote this subgroup of $G$ by $H$. Choose $k\in\{1,2,3,4\}\sm\{i,j\}$. Then
\begin{equation*}\label{eq-3.7}
H_i\cap H_j=(G_i\cap G_j)(G_i\cap G_k)\cap (G_i\cap
G_j)(G_j\cap G_k)\supseteq G_i\cap G_j.
\end{equation*}
On the other hand, $H_i\cap
H_j\subseteq G_i\cap G_j$ since $H_i\subseteq G_i$ and $H_j\subseteq G_j$.
Therefore we have
\begin{equation}\label{eq-3.8}H_i\cap H_j=G_i\cap G_j.
\end{equation}
In the case $i\not=4$,
$a_iH_i\cap \tilde{a_i}H_j=\emptyset$ by Remark 3.1, and hence
$a_iH_iH_j\cap \tilde{a_i}H_iH_j=\emptyset$ by Lemma 2.1(v).
So we have
\begin{equation}\label{eq-3.9}i\in\{1,2,3\}\Rightarrow a_iH\cap \tilde{a_i}H=\emptyset.
\end{equation}

Let $1\le i\le 3$ and $j\in\{1,2,3,4\}\sm\{i\}$. As both $a_iH_i\cap a_jH_j$ and
$\tilde a_iH_i\cap a_jH_j$ are contained in $a_iG_i\cap a_jG_j=\emptyset,$
by Lemma 2.1(v) we have
$$a_iH\cap a_jH=a_iH_iH_j\cap  a_jH_iH_j=\emptyset$$
and
$$\tilde a_iH\cap a_jH=\tilde a_iH_iH_j\cap
 a_jH_iH_j=\emptyset.$$
If $1\le i<j\le 3$, then
$$\tilde a_iH_i\cap \tilde a_jH_j\se a_iG_i\cap a_jG_j=\emptyset$$
and hence $\tilde a_iH\cap \tilde a_jH=\emptyset$.

From the above we see that the following seven cosets
$$a_1H,\ a_2H,\ a_3H,\ a_4H,\ \tilde a_1H,\ \tilde a_2H,\ \tilde a_3H$$
are pairwise disjoint. Therefore $[G:H]\ge 7$.
On the other hand,
$$[G:H]=[G:H_1H_2]=\f{[G:H_2]}{[H_1:H_1\cap
H_2]}=\f{[G:H_1][G:H_2]}{[G:G_1\cap G_2]}=\f{6q_16q_2}{3q_1q_2r_{12}}=6.$$
So we get a contradiction which ends the discussion in Case 2.

The proof of Theorem 3.1 is now complete.

\medskip

\Ack. The author wishes to thank her advisor Prof. Zhi-Wei Sun, who
proposed the conjecture, for his many helpful comments.

\vspace{5mm}

Wan-Jie Zhu: Department of Mathematics, Nanjing University, Nanjing
210093, People's Republic of China

{\it E-mail address:} {\tt zhuwanjie1@163.com}


\begin{thebibliography}{S06}

\bibitem{O} K. O'Bryant, {\it On Z.-W. Sun's disjoint congruence classes conjecture},
Integers {\bf 7}(2)(2001), \#A30, 10pp. (electronic)

\bibitem{R} J. S. Rose, {\it A Course on Group Theory},
Cambridge Univ. Press, Cambridge, 1978.

\bibitem{S01} Z. W. Sun, {\it Exact $m$-covers of groups by
cosets}, European J. Combin. {\bf 22}(2001), 415--429.

\bibitem{S06} Z. W. Sun, {\it Finite covers of groups by cosets or
subgroups}, Internat. J. Math. {\bf 17}(2006), 1047--1064.

\bibitem{Suz} M. Suzki, {\it Group Theory  I-II},
Springer, New York, 1982, 1986.



\end{thebibliography}
\end{document}